\documentclass[a4paper,12pt,reqno]{amsart}
\usepackage{amsmath,amssymb,amsthm,bm}
\usepackage{caption,graphicx}
\usepackage[english]{babel}
\usepackage{amscd}
\usepackage{amsgen}
\usepackage[final]{epsfig}
\usepackage{latexsym}
\usepackage{amsfonts}
\usepackage{url}
\usepackage{slashed}
\usepackage[all]{xy}
\usepackage{here}
\usepackage{caption,graphicx}

\usepackage{caption,graphicx}

\newtheorem{thm}{Theorem}
\newtheorem{lemm}{Lemma}
\newtheorem{prop}{Proposition}
\newtheorem{cor}{Corollary}

\newtheorem{defn}{Definition}

\begin{document}

\title[Fundamental groups, slalom curves and extremal
length]{Fundamental groups, slalom curves and \\ extremal length}

\author[Burglind J\"oricke]{Burglind J\"oricke}

\address{Humboldt-University Berlin\\ Unter den Linden 6\\ 10099 Berlin\\
Germany}

\email{joericke@googlemail.com}
\dedicatory{To the memory of my teacher and collaborator
Viktor Havin, his enthusiasm and his ability to convey a
great feeling of the beauty of mathematics}

\begin{abstract}
We define the extremal length of elements of the fundamental group of the twice
punctured complex plane and give upper and lower bounds for this invariant. The
bounds differ by a multiplicative constant. The main motivation comes from
$3$-braid invariants and their application.
\end{abstract}



\maketitle

In this paper we will describe a conformal invariant for the elements of the
fundamental group $\pi_1(\mathbb {C}\setminus \{-1,1\}, 0)$ of the twice
punctured complex plane with base point $0$ and give upper and lower bounds for
this invariant. The group $\pi_1(\mathbb {C}\setminus \{-1,1\}, 0)$ is a free
group with two generators. We choose generators $a_1$ and $a_2$ so that $a_1$
is represented by a simple closed curve $ \alpha_1$ with base point $0$ which
surrounds the point $-1$ counterclockwise such that the image of the curve
except the point $0$ is contained in the left half-plane. Respectively, a
standard representative $\alpha_2$ of the generator $a_2$ surrounds the point
$1$ counterclockwise and the image of the curve except the point $0$ is
contained in the right half-plane.

The fundamental group $\pi_1  \stackrel{def}{=} \pi_1(\mathbb
{C}\setminus \{-1,1\},
0)$ is isomorphic to the relative fundamental group $\pi_1^{tr}
\stackrel{def}{=}\pi_1(\mathbb {C}\setminus \{-1,1\},
(-1,1))$ whose elements are homotopy classes of curves in $\mathbb
{C}\setminus \{-1,1\}$ with end points on the interval $(-1,1)$.
We refer to $\pi_1^{tr} $ as fundamental group with totally real
boundary values ($tr$-boundary values for short). To establish the
isomorphism one has to use that the set $(-1,1)$ is connected and
simply connected and contains $0$. In the same way
$\pi_1$ is isomorphic to the relative fundamental group with
perpendicular bisector boundary values $\pi_1^{pb}
\stackrel{def}{=}\pi_1(\mathbb {C}\setminus \{-1,1\},
i\mathbb{R})$ whose elements are homotopy classes of curves in
$\mathbb {C}\setminus \{-1,1\}$ with end points on the imaginary axis
$i\mathbb{R}$. For an element $w \in \pi_1$ we denote by $w_{tr}$,
and by
$w_{pb}$, respectively, the elements in $\pi_1^{tr}$, and in
$\pi_1^{pb}$, respectively, corresponding to $w$.

Consider a rectangle $\mathcal{R}$ with sides parallel to the axes and with
length of the horizontal sides equal to $\textsf{b}$ and length of the vertical
sides equal to $\textsf{a}$.  Recall that according to Ahlfors's definition
\cite{A1} the extremal length $\lambda(\mathcal{R})$ of such a rectangle is
equal to $\frac{\textsf{a}}{\textsf{b}}$ and its conformal module
$m(\mathcal{R})$  equals $\frac{\textsf{b}}{\textsf{a}}$. A continuous mapping
of the rectangle $\mathcal{R}$ into $\mathbb {C}\setminus \{-1,1\}$ is said to
represent $w_{tr}$ if it has a continuous extension to the closure of
$\mathcal{R}$ which maps horizontal sides to the interval $(-1,1)$ and whose
restriction to each vertical side represents $w_{tr}$. We make the respective
convention for $w_{pb}$ instead of $w_{tr}$.

We are now in the position to define the extremal length of elements
of the relative fundamental groups.

\begin{defn}\label{def1} For an element $w$ of the fundamental group
$\pi_1(\mathbb {C}\setminus \{-1,1\},
0)$ the extremal length of $w$ with perpendicular bisector boundary
values ($pb$-boundary values for short) is defined as
\begin{align*}
& \lambda_{pb}(w) \stackrel{def}{=}
\inf \{\lambda({\mathcal{R}}):\\ & \mathcal{R}  \mbox{ admits 
a holomorphic map to}\;\mathbb {C}\setminus \{-1,1\} \mbox{ that
represents}\; w_{pb}\}.
\end{align*}
\end{defn}
An analogous definition can be given for $\lambda_{tr}$.

We will give upper and lower bounds for $\lambda_{pb}$ and
$\lambda_{tr}$
differing by a multiplicative constant. This is of independent
interest for the fundamental group of the twice punctured plane, but
the main motivation was to give estimates of conformal invariants of
braids.
Recall that a pure geometric $n$-braid with base point is a
continuous mapping of the unit interval $[0,1]$ into $n$-dimensional
configuration space  $C_n(\mathbb{C}) = \{(z_1,\ldots,z_n): z_j \neq
z_k \; \mbox{for} \; j \neq k\}$ whose values at the endpoints are
equal to a given base point in $C_n(\mathbb{C})$. More geometrically,
a pure geometric $n$-braid consists of $n$ pairwise disjoint curves
in the cylinder $[0,1] \times \mathbb{C}$, each joining a point in
the top $\{1\} \times \mathbb{C}$ of the cylinder with its copy in
the bottom so that for each curve the canonical projection to the
interval $[0,1]$ is a homeomorphism. A pure $n$-braid with base point
is an isotopy class of pure geometric $n$-braids with fixed base
point.

Consider a pure geometric $3$-braid. Associate to it a curve in
$\mathbb{C}\setminus \{-1,1\}$ as follows.
For a point $z=(z_1,z_2,z_3)\in C_3(\mathbb{C})$ we denote by $M_z$
the
M\"{o}bius
transformation that maps $z_1$ to $0$, $z_3$ to $1$ and fixes
$\infty$. Then $M_z(z_2)$ omits $0$, $1$ and $\infty$. Notice that
$z_2$ is equal to the cross ratio $(z_2,z_3;z_1,\infty)
=\frac{z_2-z_1}{z_3-z_1}\cdot
\frac{z_3-\infty}{z_2-\infty}=\frac{z_2-z_1}{z_3-z_1}$.

Let $\gamma(t)\,=\, (\gamma_1(t),\gamma_2(t), \gamma_3(t)),\, t \in
[0,1],$ be a curve in $C_3(\mathbb{C)}$. 
Associate to it the curve $\mathfrak{C}(\gamma)(t) \stackrel{def}{=}
2
\,\frac{\gamma_2(t)-\gamma_1(t)}{\gamma_3(t)-\gamma_1(t)} \,-1\,,\,
t \in
[0,1],$ in $\mathbb{C}$ 
which omits the points $-1$ and $1$. If $\gamma$ is a loop with base
point $\gamma(0)=(-1,0,1)$ then $\mathfrak{C}(\gamma)$ is a loop with
base
point $\mathfrak{C}(\gamma)(0)=0$. The homotopy class of
$\mathfrak{C}(\gamma)$ in $\mathbb{C}
\setminus \{-1,1\}$ with base point $0$ depends only on the homotopy
class of $\gamma$
in the configuration space $C_3(\mathbb{C})$ with base point
$(-1,0,1)$. We obtain a surjective
homomorphism $\mathfrak{C}_*$ from the fundamental group of
$C_3(\mathbb{C})$ with base point $(-1,0,1)$ to the fundamental group
of
$\mathbb{C}\setminus \{-1,1\}$ with base point $0$. The kernel of
$\mathfrak{C}_*$ equals $\langle \Delta_3^2 \rangle$, the subgroup
of $\mathcal{B}_3$ generated by the full twist obtained by twisting
the cylinder keeping the bottom fixed and turning the top by the
angle $2\pi$. Respective facts hold for loops in $C_3(\mathbb{C})$
with specified boundary values instead of loops with a base point.
With the natural definition of the extremal length with totally real
boundary values of a pure $3$-braid $b$  this extremal length is
equal to $\lambda_{tr}(\mathfrak{C}_*(b))$. The respective fact holds
for  perpendicular bisector boundary values. The obtained invariants
are invariants of $3$-braids rather than invariants of conjugacy
classes of $3$-braids. In particular, they are finer than a  popular
invariant of braids, the entropy. Our estimates imply estimates of
the entropy of pure $3$-braids $b$ in terms of the representing word
of the image $ \mathfrak{C}_*(b)$.  Notice that the name "totally
real" and "perpendicular bisector" is motivated by the definition in
the case of braids. Details will be given in a later paper. For an
introduction to braids see e.g. \cite{KaTu}. For more information on
the conformal module, the extremal length  and entropy of braids, or
of conjugacy classes of braids, respectively, see also \cite{Jo1} and
\cite{Jo2}.

We will now lift the elements of $\pi_1^{pb}$ to the logarithmic
covering $U_{\log}$ of $\mathbb{C} \setminus \{-1,1\}$ and identify
the lifts with homotopy slalom curves. This geometric interpretation
will suggest how to estimate the extremal length with perpendicular
bisector boundary values.

The logarithmic covering of $\mathbb{C} \setminus \{-1,1\}$ is the universal
covering of the twice punctured Riemann sphere $\mathbb{P}^1 \setminus
\{-1,1\}$ with all preimages of $\infty$ under the covering map removed.
Geometrically the universal covering of $\mathbb{P}^1 \setminus \{-1,1\}$ can
be described as follows. Take copies of $\mathbb{P}^1 \setminus (-1,1)$ labeled
by the set $\mathbb{Z}$ of integer numbers. Close up each copy by attaching two
copies of $(-1,1)$, the $+$-edge (the accumulation set of points of the upper
half-plane) and the $-$-edge (the accumulation set of points of the lower
half-plane). For each $k \in \mathbb{Z}$ we glue the $+$-edge of the $k$-th
copy to the $-$-edge of the $k+1$-st copy (using the identity mapping on
$(-1,1)$ to identify points on different edges). Denote by $U_{\log}$ the set
obtained from the described covering by removing all preimages of $\infty$.

The following proposition holds.

\begin{prop}\label{prop1} The set $U_{\log}$ is conformally
equivalent to $\mathbb{C}\setminus  i\,\mathbb{Z}$. The mapping
$f_1 \circ f_2$, $f_2(z)= \frac{e^{\pi z} -1}{e^{\pi z} +1}, \, z
\in
\mathbb{C}\setminus  i\,\mathbb{Z}$, $f_1(w)=
\frac{1}{2}(w+\frac{1}{w}),\, w \in \mathbb{C}\setminus\{0\}$, is
a covering map from $\mathbb{C}\setminus  i \,\mathbb{Z}$ to
$\mathbb{C} \setminus \{-1,1\}$.

The lift of $\alpha_1$ with initial point $\frac{- i}{2} +  i
k$ is a curve which joins $\frac{- i}{2} +  i k$ with
$\frac{- i}{2} +  i (k+1)$ and is contained in the closed
left half-plane. The only points on the imaginary axis are the
endpoints.

The lift of $\alpha_2$ with initial point $\frac{- i}{2} +  i
k$ is a curve which joins $\frac{- i}{2} +  i k$ with
$\frac{- i}{2} +  i (k-1)$ and is contained in the closed
right half-plane. The only points on the imaginary axis are the
endpoints.
\end{prop}

Figure 1 shows the curves $\alpha_1$ and $\alpha_2$ which represent the
generators of the fundamental group $\pi_1(\mathbb{C} \setminus \{-1,1\} ,0)$
and their lifts under the covering maps $f_1$ and $f_2 \circ f_1$ . For $j=1,2$
the curves $\alpha_j '$ and $\alpha_j ''$ are the two lifts of  $\alpha_j$
under the double branched covering $f_1: \mathbb{C} \setminus \{0\} \to
\mathbb{C}$ with branch points $1$ and $-1$. The curve $\tilde {\alpha} _1'$ is
the lift of $\alpha _1'$ under the mapping $f_2$ with initial point
$\frac{-i}{2}$, the curve $\tilde {\alpha} _2'$ lifts $\alpha _2'$ and has
initial point $\frac{i}{2}$.

\begin{figure}[H]
\begin{center}
\includegraphics[width=9cm]{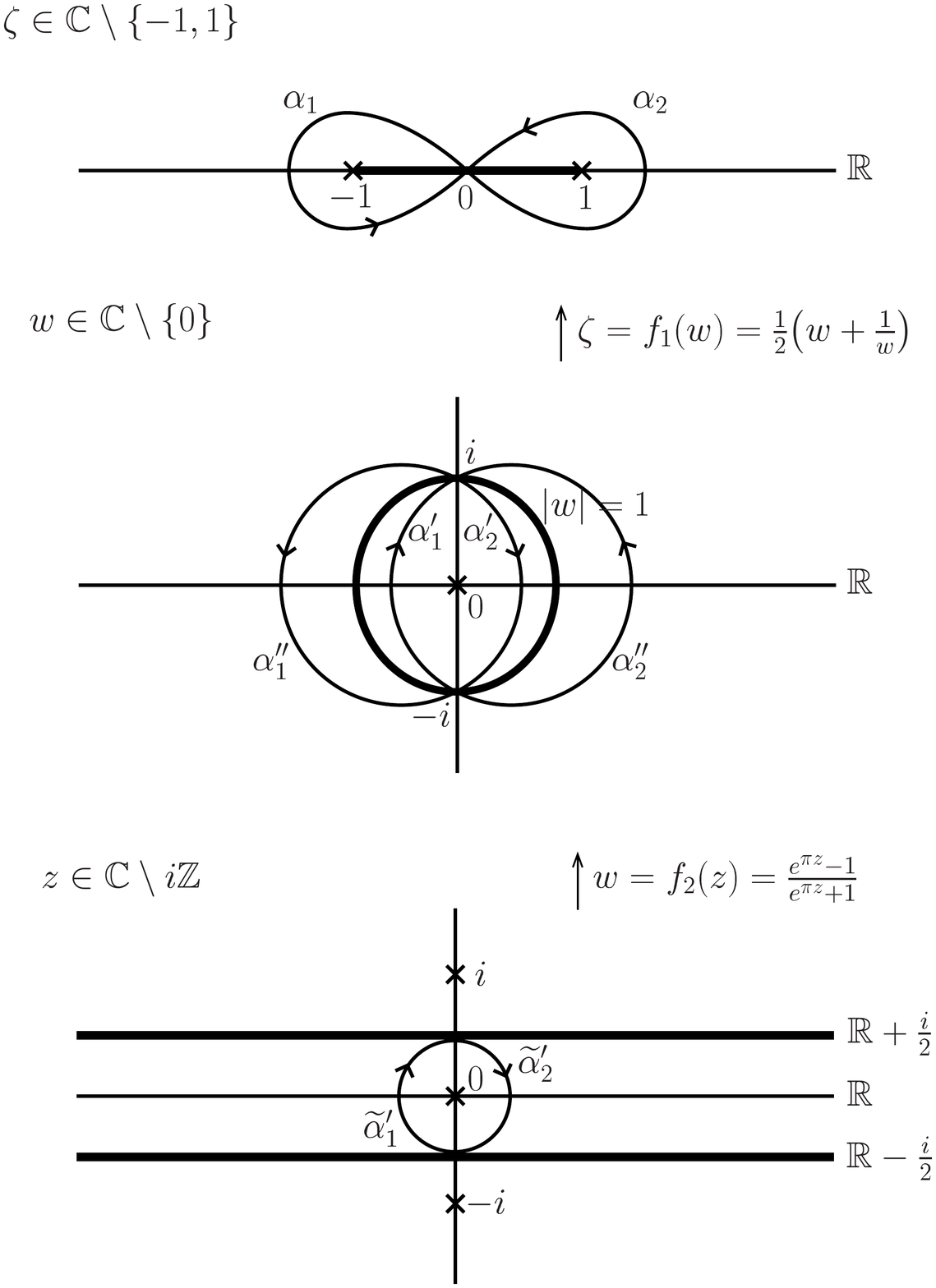}
\end{center}
\end{figure}
\centerline {Figure 1}

\bigskip

Consider the curve $\alpha_1^n,\;n\in \mathbb{Z}\setminus \{0\}$. It
runs $n$ times along the curve $\alpha_1$ if $n>0$ and $|n|$ times
along the curve which is inverse to $\alpha_1$ if $n<0$.
For each $k \in
\mathbb{Z}$ the curve $\alpha_1^n$ lifts to a curve with initial
point
 $\frac{- i}{2} +  i k$ and terminating point  $\frac{- i}{2} +  i k
 +
 i n$
which is contained in the closed left half-plane and omits the
points in $ i \mathbb{Z}$. Respectively, $\alpha_2^n,\;n\in
\mathbb{Z}\setminus \{0\},$ lifts to a curve with initial point
 $\frac{+ i}{2} +  i k$ and terminating point  $\frac{+ i}{2} +  i k
 - i n$
 which is contained in the closed right half-plane and omits the
points in $ i \mathbb{Z}$. The mentioned lifts are homotopic
through curves in $\mathbb{C} \setminus  i \mathbb{Z}$ with
endpoints on $i\,\mathbb{R} \setminus  i \mathbb{Z}$ to curves
with interior contained in the open (left, respectively, right)
half-plane. We have the following definition where we identify a curve with its image, ignoring orientation.

\begin{defn}\label{def2} A simple arc in  $\mathbb{C} \setminus  i
\mathbb{Z}$ with endpoints on different connected components of
$i\,\mathbb{R}
\setminus  i \mathbb{Z}$ is called an elementary slalom curve if
its interior (i.e. the complement of its endpoints) is contained in
one of the open half-planes $\mathbb{C}_{r}  \stackrel{def}{=}\{z \in
\mathbb{C}: \,\mbox{Re}z >0\}$ or
$\mathbb{C}_{\ell}  \stackrel{def}{=}\{z \in \mathbb{C}: \,\mbox{Re}z
<0\}$.

A curve in  $\mathbb{C} \setminus  i \mathbb{Z}$ is called an
elementary half slalom curve if one of the endpoints is contained in
a
horizontal line $\{\mbox{Im} z =k+\frac{1}{2}\}$ for an integer $k$
and the union of the curve with its mirror reflection in the line
$\{\mbox{Im} z = k+\frac{1}{2}\}$
is an elementary slalom curve.

 A slalom curve in  $\mathbb{C}
\setminus  i \mathbb{Z}$ is a curve which can be divided into a
finite number of elementary slalom curves so that consecutive
elementary slalom curves are contained in different half-planes.

A curve which is homotopic to a slalom curve in $\mathbb{C} \setminus
i
\mathbb{Z}$ through curves with endpoints in $\mathbb{R} \setminus
 i \mathbb{Z}$ is called a homotopy slalom curve.
\end{defn}

Figure 2 below shows a slalom curve which represents a lift of the
element $a_2^{-1}\,a_1^2\,
a_2^{-3}\,a_1^{-1}\,a_2^{-1}\,a_1^{-1}\,a_2\,a_1^{-1}$ with
perpendicular bisector boundary values.


\begin{figure}[H]
\begin{center}
\includegraphics[width=4cm]{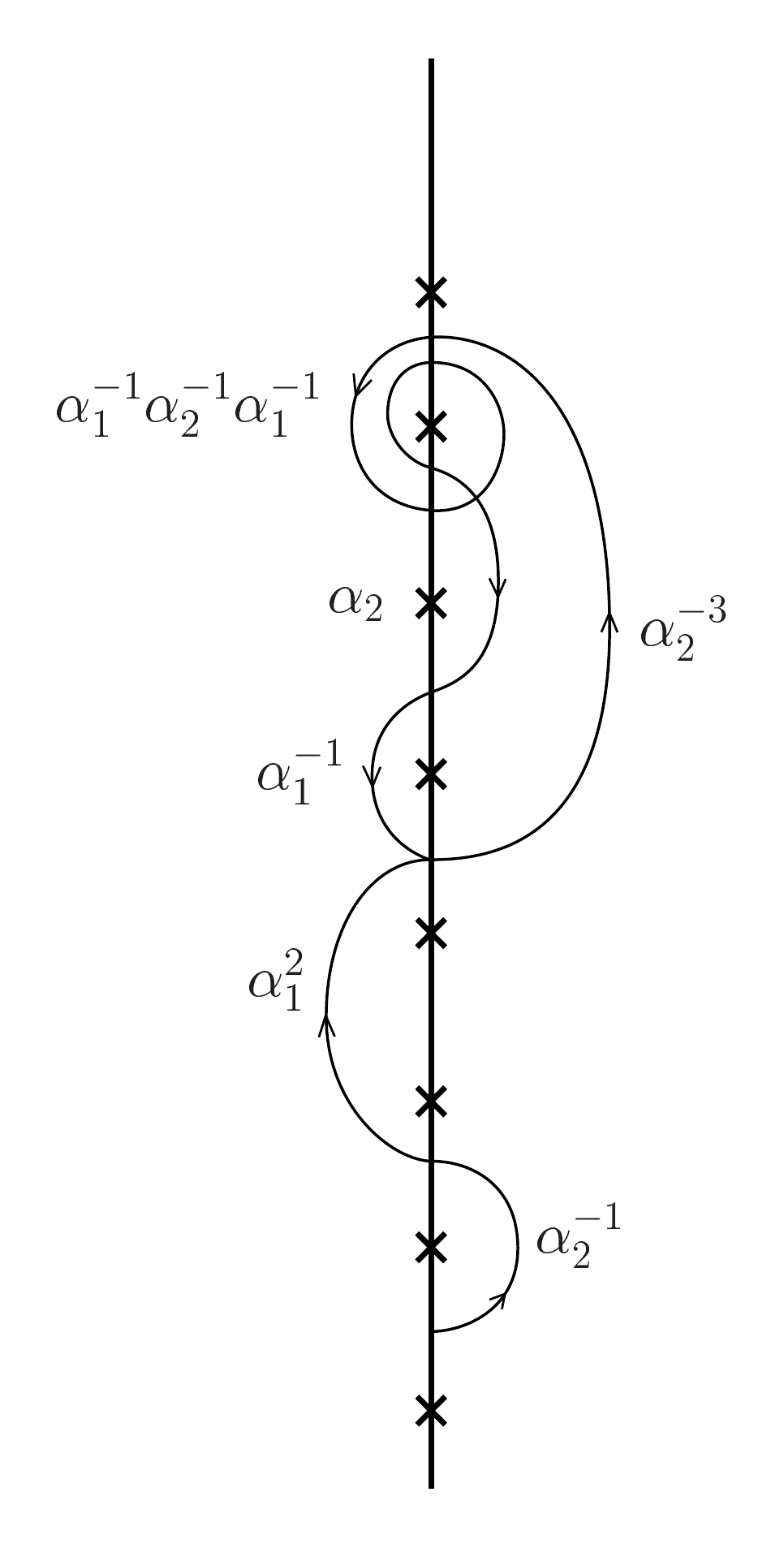}
\end{center}
\end{figure}
\centerline {Figure 2}

\bigskip

\bigskip

Elementary slalom curves and elementary half slalom curves will serve
as building blocks. Note that each curve in $\mathbb{C} \setminus  i
\mathbb{Z}$ with
endpoints in $i \mathbb{R} \setminus  i \mathbb{Z}$ is a
homotopy slalom curve or is homotopic to the identity in  $\mathbb{C}
\setminus  i \mathbb{Z}$ with
endpoints in $i \mathbb{R} \setminus  i \mathbb{Z}$. Proposition
\ref{prop1} implies that each lift of a curve in
$\mathbb{C}\setminus \{-1,1\}$ with endpoints on the imaginary axis
is of such type. The extremal length of slalom curves (more precisely
of homotopy classes of slalom curves) can be defined in the same way as
the respective object for elements of the fundamental group $ \pi_1
^{pb}$. The extremal length with perpendicular bisector boundary
values of an element of $\pi_1$ is equal to that of its lift.

Consider the extremal length of an elementary slalom curve which corresponds to
the word $a^n \in \pi_1$ where $a$ equals either $a_1$
or $a_2$. 
Without loss of generality we may assume that $a=a_1$, hence the
curve is contained in the closed left half-plane. After a translation
the endpoints of the curve are contained
in the intervals
$(-i(M+1), -iM)$ and $(iM, i(M + 1))$, respectively, with
$M=\frac{|n|-1}{2}$. For $M=0$ (thus for $|n|=1$) the extremal length
equals $0$. In this case we call the original curve a trivial elementary
slalom curve. Let $M$ be positive. The curve is represented by the
extension to the boundary of a
conformal mapping of an open rectangle $\mathcal{R}^M$ onto the left
half-plane which maps
the horizontal sides onto $[-i(M + 1), -iM]$, and $[iM, i(M + 1)]$,
respectively. Hence its extremal length is bounded from above by the
extremal length of the rectangle $\mathcal{R}^M$.

The conformal mapping of a  rectangle onto the  left half-plane
$\mathbb{C}_{\ell}$ whose extension to the boundary maps the
horizontal sides onto $[-i(M + 1), -iM]$, and $[iM, i(M + 1)]$,
respectively, is related to elliptic integrals. With a suitable
normalization of the rectangle the inverse of the mapping is
equal to the elliptic integral
\begin{align}\label{eq1}
\mathcal{F}_M(z)& = \int_0^z
\frac{d\zeta}{\sqrt{(\zeta^2-(iM)^2)(\zeta^2-(i(M+1))^2)}}
\nonumber\\
& = \frac{i}{M} \int_0^{\frac{z}{iM}}
\frac{dw}{\sqrt{(1-w^2)((1+\frac{1}{M})^2-w^2)}}, \; z \in
\mathbb{C}_{\ell}.
\end{align}

We use the branch of the square root which is positive on the
positive
real axis. The function $\mathcal{F}_M$ extends continuously to the
imaginary axis (the integral converges). The extended map maps the
closed left half-plane to a closed rectangle. The points $-i(M+1)$,
$-iM$,
$iM$ and $i(M+1)$ are mapped to the vertices of the rectangle.  It is
known and follows from formula \eqref{eq1} for the elliptic integral
that for $M \geq \frac{1}{2}$ the extremal length of the rectangle
$\mathcal{R}^M$ satisfies the inequalities
\begin{equation}\label{eq2}
c \log(1+M)\leq \lambda(\mathcal{R}^M) \leq C \log(M+1)
\end{equation}
for positive constants $c$ and $C$ not depending on $M$.

Equation \eqref{eq2} suggests the following proposition.

\begin{prop}\label{prop2} 
The extremal length $\lambda_{k,\ell}$ of an elementary slalom curve
with endpoints in the intervals $(ik, i(k + 1))$ and $(i\ell, i(\ell
+
1))$, respectively, with $|k-\ell|\geq 2$, satisfies the inequalities
\begin{equation}\label{eq3}
c' \log(1+\frac{|k-\ell|-1}{2}) \leq \lambda_{k,\ell} \leq
C'\log(1+\frac{|k-\ell|-1}{2})
\end{equation}
for positive constants $c'$ and $C'$ not depending on $k$ and $\ell$.

\end{prop}
There are explicit estimates for the constants $c'$ and $C'$.

The proof will be given elsewhere. Here we already discussed the
estimate from above. The estimate from below is more subtle. The
first difficulty is that the representing mappings for an elementary
slalom curve are not necessarily conformal mappings, they are merely
holomorphic. The second difficulty is
that the image of the rectangle is not necessarily contained in the
half plane. We can only say about the mapping that it lifts to a
holomorphic mapping into the universal covering of
$\mathbb{C}\setminus i \mathbb{Z}$ with specified boundary values.
The universal covering is a half-plane, but the horizontal sides of
the rectangle are not mapped any more into boundary intervals of the
half-plane but into some curves in the half-plane. One tool for
dealing with these difficulties is an analog of the following lemma
which is of independent interest.

\begin{lemm}\label{lemm1} Let $\mathcal{R}_1$ and $\mathcal{R}_2$ be
rectangles with sides parallel to the axes. Suppose $S_2$ is a
vertical strip bounded by the two vertical lines which are
prolongations of the vertical sides of the rectangle $\mathcal{R}_2$.
Let $f:\mathcal{R}_1 \to S_2$ be a holomorphic map whose extension to the closure maps the two
horizontal sides of $\mathcal{R}_1$ into different horizontal sides of
$\mathcal{R}_2$. Then
$$
\lambda (\mathcal{R}_1) \geq \lambda (\mathcal{R}_2) \,.
$$
Equality holds if and only if the mapping is a surjective conformal
map from $\mathcal{R}_1$ to $\mathcal{R}_2$.
\end{lemm}

The proof of the lemma is based on the Cauchy-Riemann equations.

\medskip

To estimate the extremal length of an arbitrary element of the
fundamental group $\pi_1
=\pi_1(\,\mathbb{C}\setminus
\{-1,1\}\,,\; \{0\}\,\})$  we represent the element as a word in the
generators (and identify it with the word). A word is in reduced form
(or a reduced word) if it is written as product of powers of
generators where consecutive terms correspond to different
generators. Consider the reduced word
\begin{equation}\label{eq4}
w=a_1^{n_1} \cdot
a_2^{n_2} \cdot \ldots ,\,
\end{equation}
where the $n_j$ are integers. (Here
$a_j^0  \stackrel{def}{=}\mbox{id}$, we allow $n_1=0$.) We are
interested first in the extremal length with perpendicular boundary value
conditions. 
One can show that any curve which represents this element with
perpendicular boundary values can be represented as composition of
curves $\alpha_1^{n_1}$, $\alpha_2^{n_2}$,..., which represent
$a_1^{n_1}$,  $a_2^{n_2}$,...,  with perpendicular boundary values.
Together with Theorems 2 and 4 of \cite{A1} this implies the
following estimate from below
$$
\lambda_{pb}(a_1^{n_1}
a_2^{n_2}\ldots) \geq \lambda_{pb}(a_1^{n_1})
+\lambda_{pb}(a_2^{n_2})\, + \ldots.
$$
This gives a good lower bound if all terms $a_j^{n_j}$ of the reduced word
enter with power of absolute value at least $2$. It does not give a good lower
bound, for example, for the word $(a_1 \, a_2^{-1})^n$ with $n\geq 1$,  or for
the word $a_1^{n_1} (a_2 a_1)^{n_2}  $ for integers $n_1$ and $n_2$ larger than
$1$ and $n_2$ much bigger than $n_1$. In the first example the reason is the
following. Each representing curve for $a_1 \, a_2^{-1}$ with $pb$-boundary
values can be written as composition of the following two curves: a curve
$\alpha_1$ with $pb$-boundary values on the left and $tr$-boundary values on
the right representing $a_1$, and a curve $\alpha_2^{-1}$ with $tr$-boundary
values on the left and $pb$-boundary values on the right representing
$a_2^{-1}$. The lift of each of the two curves is a non-trivial half-slalom
curve. Hence the extremal length with $pb$-boundary values of the element $(a_1
\, a_2^{-1})^n$ is proportional to $n$.

For the second example one can show that each representing curve with
$pb$-boundary values contains a piece corresponding to $(a_2
a_1)^{n_2}$ with mixed boundary values. A different choice of a lift
gives a half slalom curve which shows that the extremal length of
this piece is proportional to $\log(\frac{n_2 -1}{2})$.

The discussion suggests that the extremal length of a general element
of $\pi_1^{pb}$ can be given in terms of a syllable decomposition of
the representing reduced word.

We describe now the syllable decomposition of the word \eqref{eq4}.
\begin{itemize}
\item [(1)] Any term $a_j^{n_j}$ of the reduced word with $|n_j| \ge 2$ is a syllable.
\item [(2)] Any maximal sequence of consecutive terms of the
reduced word which have equal power equal to either $+1$ or $-1$ is a syllable.
\item [(3)] Each remaining term of the reduced word is
characterized by the following properties. It enters with power $+1$ or $-1$
and the neighbouring term on the right (if there is one) and also the
neighbouring term on the left (if there is one) has power different from that
of the given one. Each term of this type is a syllable, called a singleton.
\end{itemize}

Define the degree of a syllable $\mbox{deg}(\mbox{syllable})$ to be
the sum of the absolute values of the powers of terms entering the
syllable.


For example, the syllables of the word $a_2^{-1}\,a_1^2\,
a_2^{-3}\,a_1^{-1}\,a_2^{-1}\,a_1^{-1}\,a_2\,a_1^{-1}$ (see Figure 2) from left
to right are the singleton $a_2^{-1}$, the syllable $a_1^2$ of degree $2$, the
syllable $a_2^{-3}$ of degree $3$, the syllable $a_1^{-1}\,a_2^{-1}\,a_1^{-1}$
of degree $3$, the singleton $a_2$ and the singleton $a_1^{-1}$.

Put $\Lambda (w)  \stackrel{def}{=} \sum_{\mbox {\Small syllables
of}\;w} \log (1 \,+\, \mbox{deg (syllable)})$.

The following theorem holds.

\begin{thm}\label{thm1} There are absolute
positive constants $C_+$ and $C_-$ such that the following holds. Let
$w$ be the word representing an element of
$\pi_1= \pi_1(
\mathbb{C}\setminus\{-1,1\}\,,\;\{0\}\,)$. Then
\begin{itemize}
\item [(1)] $C_- \cdot \Lambda(w) \leq
\lambda_{tr}(w)\leq C_+ \cdot
\Lambda(w)$, except in the following cases: $w=a_1^n$ or $w=a_2^n$
for an integer $n$. In these cases $\lambda_{tr}(b)=0$.

\item [(2)]  $C_- \cdot \Lambda(w) \leq \lambda_{pb}(w) \leq C_+ \cdot
\Lambda(w)$, except in the following case: 
each term in the reduced word $w$ has the same power, which equals
either $+1$ or $-1$. In these cases $\lambda_{pb}(w)=0$.

\end{itemize}
\end{thm}

\begin{cor}\label{cor1} For an element $w \in \pi_1$
which is not one of the exceptional cases of Theorem \ref{thm1} the
two versions of the extremal length are comparable:
$$
C_1 \,\lambda_{tr}(w) \,\leq \,\lambda_{pb}(w)\,\leq \,C_2
\,\lambda_{tr}(w)
$$
for positive constants $C_1$ and $C_2$ which do not depend on $w$.
\end{cor}

\begin{cor}\label{cor2} There are positive constants $C_-'$ and
$C_+'$ such that for each element $w \in \pi_1$
which is not a singleton the estimate
$$
C_-' \cdot \Lambda(w) \leq \lambda_{tr}(w) + \lambda_{pb}(w) \leq C_+
' \cdot \Lambda(w)
$$
holds.
\end{cor}

The extremal length of elements of the fundamental group of the complex plane
with an arbitrary number of punctures will be treated in a forthcoming paper.
The case of $n$-braids with arbitrary $n$ is more subtle.

During the work on the present paper the author was supported by the SFB
"Space-Time-Matter " at Humboldt-University Berlin. She is grateful to A.
Khrabrov and  K. Mshagskiy  for the professional drawing of the figures.

\end{document}